\documentclass[12pt,a4paper,twoside]{article}

\newcommand{\pageformat}[6]{\setlength{\hoffset}{-1in}
                  \setlength{\voffset}{-1in}
                  \addtolength{\hoffset}{#5}
                            \addtolength{\voffset}{#6}
                            \setlength{\oddsidemargin}{#1}
                            \setlength{\evensidemargin}{#2}
                            \setlength{\textwidth}{\paperwidth}
                  \addtolength{\textwidth}{-\oddsidemargin}
                  \addtolength{\textwidth}{-\evensidemargin}
                  \addtolength{\textwidth}{-\marginparsep}
                  \addtolength{\textwidth}{-\marginparwidth}
                            \setlength{\topmargin}{#3}
                            \setlength{\textheight}{\paperheight}
                  \addtolength{\textheight}{-\topmargin}
                  \addtolength{\textheight}{-\headheight}
                  \addtolength{\textheight}{-\headsep}
                  \addtolength{\textheight}{-\footskip}
                  \addtolength{\textheight}{-#4}}
\pageformat{2cm}{3cm}{24mm}{24mm}{1pt}{0pt}

\usepackage{ifthen}
\newboolean{article}
    \setboolean{article}{true}
\newboolean{report}
\newboolean{book}
\newboolean{letter}
\newboolean{german}
\newboolean{italian}
\newboolean{nobaselinestretch}
\newboolean{nosectionappendix}
\newboolean{oldtoc}
\newboolean{nosectionequn}
\newboolean{notheorem}

\ifthenelse{\boolean{german}}{
    \usepackage{german}}{}

\usepackage[latin1]{inputenc}

\usepackage{amsmath}
\usepackage{amssymb}
\usepackage[mathscr]{eucal}

\ifthenelse{\boolean{notheorem}}{}{
    \usepackage{theorem}}

\ifthenelse{\boolean{nobaselinestretch}}{}{
    \renewcommand{\baselinestretch}{1.25}}

\newenvironment{env}[2]{\begin{#1}#2\end{#1}}{}
    \newcommand{\beq}[1]{\begin{env}{equation}{#1}}
    \newcommand{\beqn}[1]{\begin{env}{equation*}{#1}}
    \newcommand{\bal}[1]{\begin{env}{align}{#1}}
    \newcommand{\baln}[1]{\begin{env}{align*}{#1}}
    \newcommand{\bga}[1]{\begin{env}{gather}{#1}}
    \newcommand{\bgan}[1]{\begin{env}{gather*}{#1}}
    \newcommand{\bflal}[1]{\begin{env}{flalign}{#1}}
    \newcommand{\bflaln}[1]{\begin{env}{flalign*}{#1}}
    \newcommand{\bmu}[1]{\begin{env}{multline}{#1}}
    \newcommand{\bmun}[1]{\begin{env}{multline*}{#1}}
    \newcommand{\bsp}[1]{\begin{env}{split}{#1}}

    \newcommand{\eeq}{\end{env}}
    \newcommand{\eeqn}{\end{env}}
    \newcommand{\eal}{\end{env}}
    \newcommand{\ealn}{\end{env}}
    \newcommand{\ega}{\end{env}}
    \newcommand{\egan}{\end{env}}
    \newcommand{\eflal}{\end{env}}
    \newcommand{\eflaln}{\end{env}}
    \newcommand{\emu}{\end{env}}
    \newcommand{\emun}{\end{env}}
    \newcommand{\esp}{\end{env}}

\newcommand{\lf}{\vspace{2ex}}

\renewcommand{\bf}[1]{\textbf{#1}}
\renewcommand{\it}[1]{\textit{#1}}

\renewcommand{\sf}[1]{\textsf{#1}}

\renewcommand{\tt}[1]{\texttt{#1}}
\newcommand{\hl}[1]{\bf{\it{#1}}}

\newcommand{\mbf}[1]{\mathbf{#1}}
\newcommand{\msf}[1]{\text{\small$\sf{#1}$}}

\newcommand{\cmc}[1]{\mathcal{#1}}
\newcommand{\eus}[1]{\mathscr{#1}}
\newcommand{\euf}[1]{\mathfrak{#1}}
\newcommand{\bb}[1]{\mathbb{#1}}

\newcommand{\nbd}[1]{$#1$\nobreakdash--}
\newcommand{\ol}[1]{\overline{#1}}

\newcommand{\vk}{\varkappa}

\newcommand{\om}{\omega}

\newcommand{\abs}[1]{\left\lvert#1\right\rvert}
\newcommand{\norm}[1]{\left\lVert#1\right\rVert}

\newcommand{\bnorm}[1]{\bigl\lVert#1\bigr\rVert}

\newcommand{\snorm}[1]{\norm{\smash{#1}}}

\newcommand{\bfam}[1]{\bigl(#1\bigr)}

\newcommand{\AB}[1]{\langle#1\rangle}

\newcommand{\CB}[1]{\{#1\}}
\newcommand{\bCB}[1]{\bigl\{#1\bigr\}}
\newcommand{\BCB}[1]{\Bigl\{#1\Bigr\}}
\newcommand{\SB}[1]{[#1]}

\newcommand{\set}[2][]{
    \ifthenelse{\equal{#1}{}}{
        \CB{#2}}{
        \CB{#1~|~#2}}}
\newcommand{\bset}[2][]{
    \ifthenelse{\equal{#1}{}}{
        \bCB{#2}}{
        \bCB{#1~|~#2}}}
\newcommand{\Bset}[2][]{
    \ifthenelse{\equal{#1}{}}{
        \BCB{#2}}{
        \BCB{#1~\big|~#2}}}

\DeclareMathOperator{\ls}{\normalfont\msf{span}}
\DeclareMathOperator{\cls}{\ol{\ls}}

\DeclareMathOperator{\id}{\normalfont\msf{id}}

\newcommand{\C}{\bb{C}}

\newcommand{\bJ}{\bb{J}}

\newcommand{\N}{\bb{N}}

\newcommand{\R}{\bb{R}}

\newcommand{\cB}{\cmc{B}}

\newcommand{\sB}{\eus{B}}

\newcommand{\et}{\euf{t}}

\newcommand{\eB}{\euf{B}}

\newcommand{\eH}{\euf{H}}

\newcommand{\eK}{\euf{K}}
\newcommand{\eL}{\euf{L}}

\newcommand{\eU}{\euf{U}}

\newcommand{\U}{\mbf{1}}

\newcommand{\G}{\Gamma}

\newcommand{\I}{{I\!\!\!\;I}}

\ifthenelse{\boolean{nosectionequn}}{}{
    \numberwithin{equation}{section}
    }

\ifthenelse{\boolean{article}\or\boolean{letter}\or\boolean{nosectionequn}}{
    \setboolean{nosectionappendix}{true}}{}
\ifthenelse{\boolean{nosectionappendix}}{}{
    \renewcommand{\appendix}{
        \chapter*{\appendixname}
        \addcontentsline{toc}{chapter}{\appendixname}
        \renewcommand{\thesection}{\Alph{section}}
        \setcounter{section}{0}}}
   
\ifthenelse{\boolean{report}\or\boolean{book}}{
    }{}

\ifthenelse{\boolean{notheorem}}{}{
        \newcommand{\mnname}{Mathematical note.}
        \newcommand{\enname}{End of the note.}
        \newcommand{\definame}{Definition.}
        \newcommand{\propname}{Proposition.}
        \newcommand{\lemname}{Lemma.}
        \newcommand{\exname}{Example.}
        \newcommand{\exername}{Exercise.}
        \newcommand{\remname}{Remark.}
        \newcommand{\obname}{Observation.}
        \newcommand{\thmname}{Theorem.}
        \newcommand{\corname}{Corollary.}
        \newcommand{\proofname}{Proof.}
    \ifthenelse{\boolean{german}}{
        \renewcommand{\mnname}{Mathematische Notiz.}
        \renewcommand{\enname}{Ende der Notiz.}
        \renewcommand{\exname}{Beispiel.}
        \renewcommand{\exername}{Übung.}
        \renewcommand{\remname}{Bemerkung.}
        \renewcommand{\obname}{Beobachtung.}
        \renewcommand{\thmname}{Satz.}
        \renewcommand{\corname}{Korollar.}
        \renewcommand{\proofname}{Beweis.}}{}
    \ifthenelse{\boolean{italian}}{
        \renewcommand{\mnname}{Nota matematica.}
        \renewcommand{\enname}{Fina della nota.}
        \renewcommand{\definame}{Definizione.}
        \renewcommand{\propname}{Proposizione.}
        \renewcommand{\exname}{Esempio.}
        \renewcommand{\exername}{Esercizio.}
        \renewcommand{\remname}{Nota.}
        \renewcommand{\obname}{Osservazione.}
        \renewcommand{\thmname}{Teorema.}
        \renewcommand{\corname}{Corollario.}
        \renewcommand{\proofname}{Dimostrazione.}

       \renewcommand{\appendixname}{Appendice}

       }{}
    \theoremheaderfont{\normalfont\bfseries}
    \theoremstyle{change}
        \theorembodyfont{\rmfamily}
            \newtheorem{emp}{}[section]
                \newcommand{\bemp}[1][]{
                    \begin{emp}\hskip-\labelsep\bf{#1}\hskip\labelsep}
                \newcommand{\eemp}{\end{emp}}
\newtheorem{itemp}[emp]{}
                \newcommand{\bitemp}[1][]{
                    \begin{itemp}\hskip-\labelsep\bf{#1}\hskip\labelsep\normalfont\itshape}
                \newcommand{\eitemp}{\end{itemp}}
            \newtheorem{mn}[emp]{\mnname}
                \newcommand{\bnm}{\begin{mn}~\begin{quotation}\renewcommand{\baselinestretch}{1}\small\noindent\ignorespaces}
                \newcommand{\enm}{\end{quotation}\hfill\bf{\enname}\end{mn}}
            \newtheorem{ex}[emp]{\exname}
                \newcommand{\bex}{\begin{ex}}
                \newcommand{\eex}{\end{ex}}
            \newtheorem{exer}[emp]{\exername}
                \newcommand{\bexer}{\begin{exer}}
                \newcommand{\eexer}{\end{exer}}
            \newtheorem{defi}[emp]{\definame}
                \newcommand{\bdefi}{\begin{defi}}
                \newcommand{\edefi}{\end{defi}}
            \newtheorem{rem}[emp]{\remname}
                \newcommand{\brem}{\begin{rem}}
                \newcommand{\erem}{\end{rem}}
            \newtheorem{ob}[emp]{\obname}
                \newcommand{\bob}{\begin{ob}}
                \newcommand{\eob}{\end{ob}}
        \theorembodyfont{\normalfont\itshape}
            \newtheorem{thm}[emp]{\thmname}
                \newcommand{\bthm}{\begin{thm}}
                \newcommand{\ethm}{\end{thm}}
            \newtheorem{prop}[emp]{\propname}
                \newcommand{\bprop}{\begin{prop}}
                \newcommand{\eprop}{\end{prop}}
            \newtheorem{cor}[emp]{\corname}
                \newcommand{\bcor}{\begin{cor}}
                \newcommand{\ecor}{\end{cor}}
            \newtheorem{lem}[emp]{\lemname}
                \newcommand{\blem}{\begin{lem}}
                \newcommand{\elem}{\end{lem}}
\newenvironment{empn}[1]{\lf\noindent\bf{#1}\ignorespaces\hskip\labelsep}{\lf}
		\newcommand{\bempn}[1]{\begin{empn}{#1}}
		\newcommand{\eempn}{\end{empn}}
		\newcommand{\bitempn}[1]{\begin{empn}{#1}\normalfont\itshape}
		\newcommand{\eitempn}{\end{empn}}
                \newcommand{\bnmn}{\begin{empn}{\mnname}~\begin{quotation}\renewcommand{\baselinestretch}{1}\small\noindent\ignorespaces}
                \newcommand{\enmn}{\end{quotation}\hfill\bf{\enname}\end{empn}}
		\newcommand{\bexn}{\begin{empn}{\exname}}
		\newcommand{\eexn}{\end{empn}}
		\newcommand{\bexern}{\begin{empn}{\exername}}
		\newcommand{\eexern}{\end{empn}}
		\newcommand{\bdefin}{\begin{empn}{\definame}}
		\newcommand{\edefin}{\end{empn}}
		\newcommand{\bremn}{\begin{empn}{\remname}}
		\newcommand{\eremn}{\end{empn}}
		\newcommand{\bobn}{\begin{empn}{\obname}}
		\newcommand{\eobn}{\end{empn}}

\newcommand{\qedsymbol}{~\rule[-0.35mm]{2mm}{2mm}}
    \newcounter{proof}[emp]
    \newenvironment{Proof}[1]{
        \vspace{1ex}
        \renewcommand{\item}[1][\stepcounter{proof}(\roman{proof})]%
            {##1\hskip\labelsep}
        \noindent\textsc{#1\hskip\labelsep}}{
        \nolinebreak\qedsymbol}
    \newcommand{\proof}[1][\proofname]{
        \begin{Proof}{#1}\ignorespaces}
    \newcommand{\qed}{\end{Proof}}
    \newcommand{\noqed}{
        \renewcommand{\qedsymbol}{}
        \end{Proof}}}
    \ifthenelse{\boolean{italian}}{
        \renewcommand{\proofname}{Dimostrazione.}}{}

\usepackage[varg]{txfonts}

\begin{document}




\title{Constructing Units in Product Systems}
\author{Volkmar Liebscher and Michael Skeide\thanks{MS is supported by research fonds of the University of Molise.}}

\date{Campobasso, October 2005}

{
\renewcommand{\baselinestretch}{1}
\maketitle



\begin{abstract}\noindent
We prove a criterion that allows to construct units in product systems of correspondences with prescribed infinitesimal characterizations. This criterion summarizes proofs of known results and new applications. It also frees the hypothesis from the assumption that the units are contained in a product system of time ordered Fock modules.
\end{abstract}


}


%

\section{Introduction} \label{intro}

{\parskip0.5ex plus 0.5ex minus 0.5ex
An \it{Arveson system} is a (measurable) family $H^\otimes=\bfam{H_t}_{t\in\R_+}$ of (infinite dimensional when $t>0$ and separable) Hilbert spaces with a (measurable) associative identification $H_s\otimes H_t=H_{s+t}$. A \it{unit} is a (measurable nontrivial) section $u^\otimes=\bfam{u_t}_{t\in\R_+}$ that factors according to that identification: $u_s\otimes u_t=u_{s+t}$. Already in his trailblazing paper \cite{Arv89} Arveson considered a sort of \it{Trotter product} composing two units $u^\otimes,v^\otimes$ to give a new one $w^\otimes$, whose elements are defined as
\beq{\label{PA1}
w_t
~:=~
\lim_{n\to\infty}\bfam{u_{\frac{t}{2n}}\otimes v_{\frac{t}{2n}}}^{\otimes n}.
}\eeq
The new unit is charcacterised infinitesimally by the values of the so-called \it{covariance function} defined as $\gamma^{u,v}:=\frac{d}{dt}\big|_{t=0}\AB{u_t,v_t}$, namely,
\beq{\label{PA2}
\gamma^{x,w}
~=~
{\textstyle\frac{1}{2}}\gamma^{x,u}
+
{\textstyle\frac{1}{2}}\gamma^{x,v}
}\eeq
for every unit $x^\otimes$.

However, when we wish to apply the idea for such type of constructions of new units to product systems of Hilbert modules then we meat two obstacles. Firstly, the limit in \eqref{PA1} is only weak. In Hilbert spaces this is not a big problem, because weak and norm closures of subspaces coincide and proofs of the \it{Mazur theorem} even show how to transform a weakly convergent sequence into a norm convergent sequence. In Hilbert modules having only convergence of inner products is deadly. Secondly, the best modification of the limit in \eqref{PA1} allows to have convex combinations on the infinitesimal level in \eqref{PA2}. (For positive $\vk,\lambda$ with $\vk+\lambda=1$ consider $w_t=\lim_{n\to\infty}\bfam{u_{\frac{\vk t}{n}}\otimes v_{\frac{\lambda t}{n}}}^{\otimes n}$. Then $\gamma^{x,w}=\vk\gamma^{x,u}+\lambda\gamma^{x,v}$.) In critical applications like, for instance, the \it{Trotter product} of units in \it{spatial} product systems as defined in Skeide \cite{Ske01p}, convex combination is not enough. We need affine combinations, that is, $\vk,\lambda$ complex with $\vk+\lambda=1$.

In these notes we prove a powerful criterion (split into two lemmata) that allows to subsummarize all existing constructions of units in a general abstract scheme. In that way we avoid to have to repeat the same sort of argument in every individual case just because the formulation is not sufficiently general. The limits will be in norm. And we do not use the hypothesis (either by explicit assumption as in \cite{Ske01p} or by refering to results like in Barreto, Bhat, Liebscher and Skeide \cite{BBLS04} that establish equivalence to deep results like Christensen-Evans \cite{ChrEv79}) that the product system is contained in a subsystem of time ordered Fock modules. This makes, in fact, a couple of results proved in \cite{BBLS04} independent of \cite{ChrEv79}.
}

\section{Units and CPD-semigroups}\label{UniCPDsec}

In this section we repeat the definitions and results (mainly from \cite{BBLS04}) we need to formulate and proof the main lemma. Throughout $\cB$ denotes a unital \nbd{C^*}algebra. The concept of product systems and their relation to \nbd{E_0}semigroups makes sense also for nonunital \nbd{C^*}algebras (see Muhly, Skeide and Solel \cite{Ske02,MSS03p}), but not the concept of units; see the discussions in Skeide \cite{Ske04p,Ske05p1}.

A \hl{correspondence} over $\cB$ is Hilbert \nbd{\cB}module with a nondegenerate (that is, unital in our case) left action of $\cB$. A \hl{product system} of Hilbert modules (or of correspondences) is a family $E^\odot=\bfam{E_t}_{t\in\R_+}$ of correspondence $E_t$ over a (unital) \nbd{C^*}algebra $\cB$ with an associative identification
\beqn{
E_s\odot E_t
~=~
E_{s+t}
}\eeqn
and $E_0=\cB$ (where $\cB$ is the \hl{trivial} correspondence over $\cB$ with the natural bimodule operation and inner product $\AB{b,b'}=b^*b'$) and $E_0\odot E_t=E_t=E_t\odot E_0$ via $b\odot x_t=bx_t$, $x_t\odot b=x_tb$. We do not pose any continuity or measurability condition. (Those will be encoded into continuity properties of the units. For a definition of \it{continuous} product system see Skeide \cite{Ske03b}.)

A \hl{unit} in a product system $E^\odot$ is a family $\xi^\odot=\bfam{\xi_t}_{t\in\R_+}$ of elements $\xi_t\in E_t$ that factors as
\beqn{
\xi_s\odot\xi_t
~=~
\xi_{s+t}
}\eeqn
and $\xi_0=\U\in\cB=E_0$. The \hl{trivial} unit is $0$ for all $t>0$. For every $t>0$ we define the $\bJ_t=\CB{(t_n,\ldots,t_1)\colon\abs{\et}=t~(n\in\N,t_i>0)}$, where the \hl{length} of a tuple $\et=(t_n,\ldots,t_1)$ is $\abs{\et}:=t_n+\ldots+t_1$, while the \hl{norm} of $\et$ is $\norm{\et}:=\max(t_n,\ldots,t_1)$. We note that $\bJ_t$ becomes a lattice when we consider $\et$ as the interval partition $\bfam{\sum_{i=1}^nt_i,\ldots,\sum_{i=1}^1t_i}$; see \cite[Proposition 4.1]{BhSk00}. Clearly, for an arbitrary set $S$ of units the spaces
\beqn{
E^S_t
~:=~
\cls\bCB{b_n\xi^n_{t_n}\odot\ldots\odot b_1\xi^1_{t_1}b_0\colon n\in\N,\et\in\bJ_t,{\xi^i}^\odot\in S,b_i\in\cB}
}\eeqn
and $E^S_0=\cB$ form a product subsystem of $E^\odot$, the subsystem \hl{generated} by $S$.

The covariance function of an Arveson system, actually, is the generator of a semigroup (under pointwise multiplication) of \nbd{\C}valued \it{positive definite kernels} on the set of units and, therefore, a \it{conditionally positive definite kernel}. For product systems of Hilbert modules the concept of \nbd{\cB}valued positive definite kernels makes sense but is no longer sufficient if we wish to speak about semigroups. Our kernels assume values in $\sB(\cB)$, the bounded linear mappings on $\cB$. A kernel $\eK\colon S\times S\rightarrow\sB(\cB)$ on a set $S$ is \hl{completely positive definite}, if
\beq{\label{cpd}
\sum_{i,j=1}^nb_i^*\eK^{\sigma_i,\sigma_j}(a_i^*a_j)b_j
~\ge~
0
\text{~~~for all~~~}
n\in\N;a_i,b_i\in\cB,\sigma_i\in S.
}\eeq
$\eK$ is \hl{conditionally completely positve definite}, if \eqref{cpd} holds at least if the condition $\sum_{i=1}^na_ib_i=0$ is satisfied. Thinking of a number $z\in\C$ as mapping $w\mapsto zw$ in $\sB(\C)$, the notions of positive definite and conditionally positive definite \nbd{\C}valued kernels are included. A family $\eK=\bfam{\eK_t}_{t\in\R_+}$ of completely positive definite kernels on $S$ is a \hl{CPD-semigroup} if for all $\sigma,\sigma'\in S$ the mappings $\eK^{\sigma,\sigma'}_t$ form a semigroup (with identity) in $\sB(\cB)$. The CPD-semigroup $\eK$ is \hl{uniformly continuous}, if every semigroup $\eK^{\sigma,\sigma'}$ is uniformly continuous. By \cite[Theorem 3.4.7]{BBLS04} the formula $\eK_t=e^{t\eL}$ establishes a one-to-one correspondence between uniformly continuous CPD-semigroups $\eK$ and conditionally completely positive definite kernels $\eL$. This unifies and generalizes the well-known results on the generators of semigroups of positive definite kernels and of CP-semigroups.

Every subset $S$ of units in a product system $E^\odot$ gives rise to a CPD-semigroup $\eU$ on $S$ defined by setting
\beq{\label{uniCPD}
\eU^{\xi,\xi'}_t
~=~
\AB{\xi_t,\bullet\xi'_t}.
}\eeq
We say the set $S$ of units is \hl{continuous}, if $\eU$ is uniformly continuous. A product system is \hl{type I}, if it is generated by a continuous subset of units. Recall that type I product systems need not be \it{spatial} in the sense of \cite{Ske01p} (generalizing the definition of Powers \cite{Pow88}). In fact, it is a major achievment of Section \ref{critSEC} that the results hold for general type I systems. Only type I product systems of von Neumann modules are spatial automatically; see \cite{BBLS04}.

Finally, by \cite[Theorem 4.3.5]{BBLS04} every (uniformly continuous) CPD-semigroup $\eK$ arises as the CPD-semigroup associated by \eqref{uniCPD} with a (continuous) set of units in a product system. (This is a generalization of Bhat and Skeide \cite[Theorem 4.8]{BhSk00} for CP-semigroups.) The unique minimal version among these product systems we call the \hl{GNS-system} of $\eK$.

\section{The criterion}\label{critSEC}

Suppose $t\mapsto y_t\in E_t$ is differentiable at $t=0$ in a suitable sense. It is our goal to find a criterion to check when $y_{t_n}\odot\ldots\odot y_{t_1}$ converges over the net $\bJ_t$ to $\zeta_t$ (necessarily forming a unit). In a first step (Lemma \ref{Lemma1}) we show (by rather algebraic means like GNS-construction) that there exists a product system, possibly bigger than the original one, that contains a unit $\zeta^\odot$ which fits suitably the infinitesimal characterization of $y_t$. Then (Lemma \ref{Lemma2}) we show that convergence happens, if and only if the new product system matches the old one, and give an applicable criterion for that it happens.

\blem\label{Lemma1}
Let $E^\odot$ be a product system that is generated by a continuous subset $S$ of units. Suppose there is a function $t\mapsto y_t\in E_t$ and there are mappings $K$, $K_\xi$ $(\xi^\odot\in S)$ in $\sB(\cB)$ such that
\bal{\label{y1}
\AB{y_t,\bullet y_t}
&
~=~
\id_\cB+tK+O(t^2),
\\\label{y2}
\AB{y_t,\bullet\xi_t}
&
~=~
\id_\cB+tK_\xi+O(t^2).
}\ealn
Then there exists a product system $F^\odot$ that contains $E^\odot$ and a unit $\zeta^\odot$ such that $S\cup\CB{\zeta^\odot}$ is continuous and
\beqn{
\eL^{\zeta,\zeta}
~=~
K
\text{~~~and~~~}
\eL^{\zeta,\xi}
~=~
K_\xi.
}\eeqn
\elem

\proof
We are done, if we show that the kernel $\eL$ on $S$ extended to $S\cup\CB{\zeta^\odot}$ in the way stated in the lemma (and $\eL^{\xi,\zeta}=*\circ K_\xi\circ*$) is still conditionally completely positive definite. In this case, $e^{t\eL}$ is a uniformly continuous CPD-semigroup on $S\cup\CB{\zeta^\odot}$ and for the product system $F^\odot$ we may choose the GNS-system of that CPD-semigroup.

We define the family of completely positive definite kernels $\eK_t$ on $S\cup\CB{\zeta^\odot}$ (think of $S\cup\CB{\zeta^\odot}$ as disjoint union), by setting
\baln{
\eK^{\xi,\xi'}_t
&
~=~
\AB{\xi_t,\bullet\xi'_t}
\text{~~~for~~~}
\xi^\odot,{\xi'}^\odot\in S,
\\
\eK^{\zeta,\xi'}_t
&
~=~
\AB{y_t,\bullet\xi'_t},
\\
\eK^{\xi,\zeta}_t
&
~=~
\AB{\xi_t,\bullet y_t}
\\
\eK^{\zeta,\zeta}_t,
&
~=~
\AB{y_t,\bullet y_t}.
}\ealn
These kernels are CPD but need not form a semigroup. Nevertheless, as for CPD-semigroups one easily shows that the kernel $\eL^{\xi,\xi'}=\lim_{t\to0}\frac{\eK^{\xi,\xi'}_t-\id_\cB}{t}$, as a limit of conditionally completely positive definite kernels, is a conditionally completely positive definite kernel, too.\qed

\brem
Of course, $F^\odot$ is unique, if we require that $S\cup\CB{\zeta^\odot}$ is generating.
\erem

In Lemma \ref{Lemma2} we will show that the elements of the unit $\zeta^\odot$ can be obtained as a norm limit of
\beq{\label{ytdef}
y_\et
~:=~
y_{t_n}\odot\ldots\odot y_{t_1}
}\eeq
over the net $\bJ_t$, if and only if $\zeta^\odot\in E^\odot$, and we will provide an easily applicable necessary and sufficient criterion. But, first we draw some consequences from Lemma \ref{Lemma1} that are independent.

\bprop\label{prepprop}
Under the hypothesis of Lemma \ref{Lemma1}:
\begin{enumerate}
\item\label{P1}
$\lim\limits_{\et\in\bJ_t}\AB{y_\et,\bullet y_\et}~=~\AB{\zeta_t,\bullet\zeta_t}$.

\item\label{P2}
$\lim\limits_{\et\in\bJ_t}\AB{\xi_t,\bullet y_\et}~=~\AB{\xi_t,\bullet\zeta_t}$ for all $\xi^\odot\in S$.
\end{enumerate}
All limits are in $\norm{\et}\to0$ uniformly in $t\in\SB{0,T}$ for every $T\in\R_+$.
\eprop

\proof
By assumption there is a constant $M>0$ such that $\norm{\AB{y_t,\bullet y_t}-\id_\cB-tK}\le t^2M$ for all sufficiently small $t$. Therefore, $\norm{\AB{y_t,\bullet y_t}}\le 1+t\norm{K}+t^2M\le e^{t\max(\norm{K},M)}$, and further
\beqn{
\bnorm{\AB{y_\et,\bullet y_\et}}
~=~
\norm{\AB{y_{t_1},\bullet y_{t_1}}\circ\ldots\circ\AB{y_{t_n},\bullet y_{t_n}}}
~\le~
e^{t_1\max(\norm{K},M)}\ldots e^{t_n\max(\norm{K},M)}
~=~
e^{t\max(\norm{K},M)}
}\eeqn
for all $t$ and all $\et$ with $\norm{\et}$ sufficiently small. In other words, the net $\bfam{y_\et}_{\et\in\bJ_t}$ is eventually bounded.

Let us write $Y_t:=\AB{y_t,\bullet y_t}$ and $Z_t:=\AB{\zeta_t,\bullet\zeta_t}$. We compute
\bmun{
\AB{y_\et,\bullet y_\et}-\AB{\zeta_t,\bullet\zeta_t}
~=~
Y_{t_1}\circ\ldots\circ Y_{t_n}-Z_{t_1}\circ\ldots\circ Z_{t_n}
\\
~=~
\sum_{k=1}^nY_{t_1}\circ\ldots\circ Y_{t_{k-1}}\circ(Y_{t_k}-Z_{t_k})\circ Z_{t_{k+1}}\circ\ldots\circ Z_{t_n}.
}\emun
We have $\snorm{Y_{t_1}\circ\ldots\circ Y_{t_{k-1}}}\le e^{(t_1+\ldots+t_{k-1})\max(\norm{K},M)}$, $\snorm{Z_{t_{k+1}}\circ\ldots\circ Z_{t_n}}\le e^{(t_{k+1}+\ldots+t_n)\norm{K}}$ and $\snorm{Y_{t_k}-Z_{t_k}}$ $\le t_k^2(M+\norm{K}^2e^{t_k\norm{K}})$. Altogether, there is a constant $M'$ (not depending on $t\in\SB{0,T}$) such that
\bal{\notag
\norm{\AB{y_\et,\bullet y_t}-\AB{\zeta_t,\bullet\zeta_t}}
~\le~
&
e^{t\max(\snorm{K},M)}M'\sum_{k=1}^nt_k^2
\\\label{basest}
&
~\le~
\norm{\et}e^{t\max(\snorm{K},M)}M'\sum_{k=1}^nt_k
~=~
\norm{\et}te^{t\max(\snorm{K},M)}M'.
}\eal
This shows \ref{P1}.

Similarly, we compute
\bmun{
\AB{y_\et,\bullet\xi_t}-\AB{\zeta_t,\bullet\xi_t}
\\
~=~
\sum_{k=1}^n\AB{y_{t_1},\bullet\xi_{t_1}}\circ\ldots\circ\AB{y_{t_{k-1}},\bullet\xi_{t_{k-1}}}\circ\AB{y_{t_k}-\zeta_{t_k},\bullet\xi_{t_k}}\circ\AB{\zeta_{t_{k+1}},\bullet\xi_{t_{k+1}}}\circ\ldots\circ\AB{\zeta_{t_n},\bullet\xi_{t_n}}.
}\emun
By hypothesis for every $\xi^\odot\in S$ there is a constant $M_\xi>0$ such that $\norm{\AB{y_t-\zeta_t,\bullet\xi_t}}\le t^2M_\xi$ for all sufficiently small $t$. And, of course, $\snorm{\AB{y_\et,\bullet\xi_t}}\le\snorm{y_\et}\snorm{\xi_t}$. By an estimate very similar to \eqref{basest} we show also \ref{P2}.\qed

\blem\label{Lemma2}
For the product system $F^\odot$ in Lemma \ref{Lemma1} the following conditions are equivalent:
\begin{enumerate}
\item\label{1}
$\lim\limits_{\et\in\bJ_t}y_{t_n}\odot\ldots\odot y_{t_1}=\zeta_t$ for all $t\in\R_+$.

\item\label{2}
$\zeta^\odot\in E^\odot\subset F^\odot$, that is, if $F^\odot$ is minimal, then $F^\odot=E^\odot$.

\item\label{3}
$\lim\limits_{\et\in\bJ_t}\AB{\zeta_t,y_{t_n}\odot\ldots\odot y_{t_1}}=\AB{\zeta_t,\zeta_t}$ for all $t\in\R_+$.
\end{enumerate}
\elem

\proof
$E_t$ is complete. So, if $\zeta_t$ is the norm limit of elements in $E_t$, then $\zeta_t\in E_t$. Therefore, \ref{1} implies \ref{2}.

As elements of the form $b_n\xi^n_{t_n}\odot\ldots\odot b_1\xi^1_{t_1}b_0$ span $E_t$ and the $y_\et$ are bounded for small $\norm{\et}$, this implies that $\lim_{\et\in\bJ_t}\AB{y_\et,\bullet x}=\AB{\zeta_t,\bullet x}$ for all $x\in E_t$. Now suppose that $\zeta^\odot\in E^\odot$, that is, $\zeta_t\in E_t$ for all $t$. Therefore, \ref{2} implies \ref{3}.

By Proposition \ref{prepprop} we have $\lim_{\et\in\bJ_t}\AB{y_\et,\bullet y_t}=\AB{\zeta_t,\bullet\zeta_t}$. Therefore, if \ref{3} holds, then we find
\beq{\label{Bwnlim}
\AB{y_\et-\zeta_t,y_\et-\zeta_t}
~=~
\AB{y_\et,y_\et}-\AB{y_\et,\zeta_t}-\AB{\zeta_t,y_\et}+\AB{\zeta_t,\zeta_t}
~\longrightarrow~
0.
}\eeq
Therefore, \ref{3} implies \ref{1}.\qed

\bob\label{sequob}
If also in \ref{3} the convergence is $O(t^2)$, then all estimates in \eqref{Bwnlim} are in $\norm{\et}$ (uniformly in $t\in\SB{0,T}$ for all $T\in\R_+$). Therefore, in this case we may pass to sequences $\bfam{\et^m}_{m\in\N}$ in $\bJ_t$ with $\norm{\et^m}\to0$. The most comon example is $\et^m=(t^m_m,\ldots,t^m_1)$ with $t^m_k=\frac{1}{m}$.
\eob

\brem
The criterion in Lemma \ref{Lemma2} corresponds to the well-known method in Hilbert spaces to conclude from weak convergence and convergence of norms to convergence in norm. It cannot be applied without constructing first in Lemma \ref{Lemma1} a space that is sufficiently big to contain a candidate for the limit. Of course, we would like to give a one-step criterion allowing to check immediately norm convergence of $y_\et$ just by looking at inner products of the $y_t$ with elements in $E_t$. The failure to be able to do so underlines once more the importance of the possibility to examine properties of a CPD-semigroup in terms of its GNS-systems. Also in \cite[Theorem 4.4.12]{BBLS04} we proved an intrinsic result about CPD-semigroups by passing to the GNS-system of the CPD-semigroup.
\erem

\section{Applications}\label{appSEC}

\bemp[A counter example.~]
We consider Arveson's Trotter product given by the limit in \eqref{PA1}. Specifically, as product system we consider $H_t=\G\bfam{L^2(\SB{0,t})}$ and $u^\otimes$ with $u_t=\om_t$, the vacuum, and $v^\otimes$ with $v_t=\psi(\I_{\SB{0,t}})$, the exponential vector to the indicator function of the interval $\SB{0,t}$. By Parthasarathy and Sunder \cite{PaSu98} or Skeide \cite{Ske00a}, these two units generate the whole product system.

As indicated in Observation \ref{sequob}, to treat only the case of sequences is sufficient. For the section $y$ we read from \eqref{PA1} that $y_t=u_{\frac{t}{2}}\otimes v_{\frac{t}{2}}$ and it is clear that the assumptions from Lemma \ref{Lemma1} are fulfilled. We easily compute
\beqn{
\lim_{n\to\infty}\AB{y_{\frac{t}{n}}^{\otimes n},y_{\frac{t}{n}}^{\otimes n}}
~=~
e^{\frac{t}{2}}
\text{~~~~~~while~~~~~~}
\lim_{n\to\infty}\AB{y_{\frac{t}{n}}^{\otimes n},e^{t\alpha}\psi(c\I_{\SB{0,t}})}
~=~
e^{t(\alpha+\frac{c}{2})}
}\eeqn
for every other unit $\bfam{e^{t\alpha}\psi(c\I_{\SB{0,t}})}_{t\in\R_+}$. We easily check that the product system $H^\otimes$ contains the unit $w^\otimes=\bfam{\psi(\frac{1}{2}\I_{\SB{0,t}})}_{t\in\R_+}$ such that
\beqn{
\AB{w_t,e^{t\alpha}\psi(c\I_{\SB{0,t}})}
~=~
e^{t(\alpha+\frac{c}{2})}
}\eeqn
to which, therefore, the limit in \eqref{PA1} converges weakly. However, $\AB{w_t,w_t}=e^{\frac{t}{4}}$ is strictly smaller than the limit $e^{\frac{t}{2}}$ of the norm square of $y_{\frac{t}{n}}^{\otimes n}$, so that the limit is not a norm limit. In fact, it is easy to check that the (minimal) product system from Lemma \ref{Lemma1} is $F^\otimes$ with $F_t=\G\bfam{L^2(\SB{0,t},\C\oplus\C)}$ that contains $H_t$ , as subsystem $\G\bfam{L^2(\SB{0,t},\C\oplus 0)}$ while the unit $\zeta^\otimes$ is given by $\zeta_t=\psi\bfam{(\frac{1}{2}\oplus\frac{1}{2})\I_{\SB{0,t}}}$. Clearly, the criterion Lemma \ref{Lemma2}\eqref{3} is violated.
\eemp

\bemp[Examples from \cite{BBLS04} and \cite{Ske01p} (now without embedding into Fock modules).]
We discuss two constructions of units that have been proved in \cite{BBLS04} explicitly assuming units in a time ordered product system and in \cite{Ske01p} by first constructing an embedding into a time ordered product system. Here we give a proof based on Lemmata \ref{Lemma1} and \ref{Lemma2} without any reference to a time ordered product system. (In fact, it is Lemma \ref{Lemma1} that gives the construction of a type I product system that contains a suitable unit by rather algebraic means and it is Lemma \ref{Lemma2} that helps to find out whether this unit is contained in the original system.)

For the first construction we study immediately a multi summand version, instead of the two summand versions considered inside a time ordered product system in \cite{BBLS04,Ske01p}. Suppose that ${\xi^\ell}^\odot$ $(\ell=1,\ldots,k)$ are units in a continuous generating subset $S$ of units of a product system $E^\odot$. Let $\vk_\ell$ be complex numbers such that $\vk_1+\ldots+\vk_k=1$. Put $y_t=\vk_1\xi^1_t+\ldots+\vk_k\xi^k_t$. Then the section $y$ fulfills the assumptions of Lemma \ref{Lemma1} with
\beqn{
K
~=~
\sum_{\ell,\ell'=1}^k\bar{\vk}_\ell\vk_{\ell'}\eL^{\xi^\ell,\xi^{\ell'}}
\text{~~~~~~and~~~~~~}
K_\xi
~=~
\sum_{\ell=1}^k\bar{\vk}_\ell\eL^{\xi^\ell,\xi}.
}\eeqn
Therefore,
\beqn{
\AB{\zeta_t,y_t}
~=~
\sum_{\ell'=1}^k\vk_{\ell'}\AB{\zeta_t,\xi^{\ell'}_t}
~=~
\AB{\zeta_t,\zeta_t}+O(t^2).
}\eeqn
From this it follows like in \eqref{basest} that the criterion in Lemma \ref{Lemma2}\eqref{3} is fulfilled.

The other construction allows to ``normalize'' a given (continuous) unit $\xi^\odot$ suitably within the product subsystem generated by $\xi^\odot$. Let $\beta\in\cB$ and put $y_t=\xi_te^{t\beta}$. It follows that $y$ fulfills the assumptions of Lemma \ref{Lemma1} with
\beqn{
K(b)
~=~
\eL^{\xi,\xi}(b)+\beta^*b+b\beta
\text{~~~~~~and~~~~~~}
K_{\xi'}(b)
~=~
\eL^{\xi,\xi'}(b)+\beta^*b.
}\eeqn
Also here one checks easily that Lemma \ref{Lemma2}\eqref{3} holds. Choosing $\beta=\frac{\eL^{\xi,\xi}(\U)}{2}+ih$ ($h\in\cB$ selfadjoint but otherwise arbitrary), then the unit $\zeta^\odot$ we obtain in that way determins a unital CP-semigroup $\AB{\zeta_t,\bullet\zeta_t}$ that has a generator with the same CP-part as $\eL^{\xi,\xi}$. Obviously we obtain the same unit $\zeta^\odot$, if we start with $y_t=e^{t\beta}\xi_t$.

We discuss a further construction, so far not yet considered elsewhere. Suppose that ${\xi^\ell}^\odot$ $(\ell=0,\ldots,k)$ are a continuous set of units and choose $a_\ell,b_\ell\in\cB$ $(\ell=1,\ldots,k)$ such that $\sum_{\ell=1}^ka_\ell b_\ell=0$. We put
\beqn{
y_t
~=~
\xi^0_t+\sum_{\ell=1}^ka_\ell\xi^\ell_tb_\ell.
}\eeqn
Then $y_\et$ converges in norm to the elements $\zeta_t$ of a unit that fulfills
\beqn{
\eL^{\xi,\zeta}(b)
~=~
\sum_{\ell=1}^kb^*_\ell\eL^{\xi^\ell,\xi}(a_\ell^*b)
\text{~~~~~~and~~~~~~}
\eL^{\zeta,\zeta}(b)
~=~
\sum_{\ell,\ell'=1}^kb^*_\ell\eL^{\xi^\ell,\xi^{\ell'}}(a_\ell^*ba_{\ell'})b_{\ell'}.
}\eeqn
This allows to modify the conditionally positive definite kernel $\eL$ rather arbitrarilly by a completely positive definite kernel associated canonically with every uniformly continuous CPD-semigroup. We refer to \cite[Section 5.3]{BBLS04}, in particular Theorem 5.3.2, for details.
\eemp

\bemp[Quantum Lévy processes.]
Every quantum Lévy process (see Schürmann \cite{MSchue93}) possesses a realization on a symmetric (or time ordered) Fock space as the solution of a quantum stochastic differential equation. Only recently Franz, Schürmann and Skeide (in preparation) have shown that the vacuum vector of the Fock space is cyclic for the minimal version of the process. Here is not the place to repeat all the somewhat heavy definitions. We refer the reader to the lecture notes of Franz \cite{Fra03p} where also a sketch of a proof can be found. The proof uses exactly the techniques employed in Lemmata \ref{Lemma1} and \ref{Lemma2} to construct a sufficiently large subset of units by applying the process to the vacuum.
\eemp

\setlength{\baselineskip}{2.5ex}


\newcommand{\Swap}[2]{#2#1}\newcommand{\Sort}[1]{}
\providecommand{\bysame}{\leavevmode\hbox to3em{\hrulefill}\thinspace}
\providecommand{\MR}{\relax\ifhmode\unskip\space\fi MR }
\providecommand{\MRhref}[2]{%
  \href{http://www.ams.org/mathscinet-getitem?mr=#1}{#2}
}
\providecommand{\href}[2]{#2}

\noindent
Volkmar Liebscher: {\small\itshape Institut für Mathematik und Informatik}, {\small\itshape Ernst-Moritz-Arndt-Universität Greifswald}, {\small\itshape 17487 Greifswald, Germany},
{\small\itshape E-mail: \tt{volkmar.liebscher@uni-greifswald.de}},
\\
{\small{\itshape Homepage: \tt{http://www.math-inf.uni-greifswald.de/biomathematik/liebscher}}}

\lf\noindent
Michael Skeide: {\small\itshape Dipartimento S.E.G.e S.}, {\small\itshape Università degli Studi del Molise}, {\small\itshape Via de Sanctis}, {\small\itshape 86100 Campobasso, Italy}, {\small{\itshape E-mail: \tt{skeide@math.tu-cottbus.de}}},
\\
{\small{\itshape Homepage: \tt{http://www.math.tu-cottbus.de/INSTITUT/lswas/\_skeide.html}}}


\end{document}